\documentclass{amsart}
\usepackage[utf8]{inputenc} 
\usepackage[T1]{fontenc} 
\usepackage{graphicx} 
\usepackage{subcaption}
\usepackage{amsmath}
\usepackage{amssymb}
\usepackage{amsfonts}
\usepackage{amsthm}
\usepackage{mathrsfs}
\usepackage[all]{xy}
\usepackage{graphicx}
\usepackage{color}
\usepackage{cite}
\usepackage{url}
\usepackage{indentfirst}
\usepackage[labelfont=bf,labelsep=period,justification=raggedright]{caption}
\usepackage[english]{babel}
\usepackage[utf8]{inputenc}
\usepackage{hyperref}
\usepackage[colorinlistoftodos]{todonotes}
\usepackage{tkz-fct}
\usepackage{tikz}
\usetikzlibrary{calc}
\usepackage[enableskew]{youngtab}
\usepackage{pstricks}
\usepackage{multicol}
\usepackage{graphicx} 
\usepackage{fancyhdr}
\theoremstyle{plain}
\newtheorem{thm}{Theorem}

\newtheorem{cor}[thm]{Corollary}

\theoremstyle{definition}

\theoremstyle{remark}
\newtheorem{rem}[thm]{Remark}

\numberwithin{equation}{section}
\numberwithin{thm}{section}

\title{A Unified Transformation Formula for Ramanujan's Theta Function}
\author{Mahipal Gurram}
\begin{document}
\begin{abstract}
In this paper, we derive a unified generalization of Ramanujan's transformation identities for the theta function $f(a,b)$, originally appearing in Ramanujan's Notebooks, Parts~III and IV. Using an approach based on residue-class dissections and modular substitutions, we obtain a closed transformation formula for $f(\zeta a, \zeta b)$, where $\zeta$ is a primitive root of unity $m$. As special cases, we recover and systematically prove Ramanujan's classical results for $m=2,3,$ and $4$, including even odd dissections, cubic transformation and the compact quartic form involving complex coefficients.
\end{abstract}
\maketitle
\section{Introduction}
The study of theta functions has its origins in the classical work of Jacobi, who introduced the celebrated Jacobi triple product identity\cite{Jacobi1829}:
\begin{equation}\label{eq:jacobi}
f(a,b) = \sum_{n=-\infty}^{\infty} a^{n(n+1)/2} b^{n(n-1)/2}
 = (-a;ab)_\infty\,(-b;ab)_\infty\,(ab;ab)_\infty,
\end{equation}
where $(x;q)_\infty = \prod_{k=0}^{\infty}(1 - x q^k)$ denotes the $q$-Pochhammer symbol.
The series $f(a,b)$ is now recognized as a two-variable generalization of the Jacobi theta function, and serves as a generating function for numerous $q$-series and partition identities.
Ramanujan, in his Notebooks (as edited by Berndt \cite{BerndtIII,BerndtIV}), introduced several remarkable transformation identities satisfied by $f(a,b)$. These transformations express $f(\zeta a, \zeta b)$ in terms of combinations of theta functions with rescaled arguments, where $\zeta$ is a root of unity. For instance, at $m=2$, Ramanujan obtained the decomposition of $f(a,b)$ into even and odd parts:
\begin{align*}
S_0 &= f(a^3b, ab^3) = \tfrac{1}{2}\left[f(a,b) + f(-a,-b)\right],\\
S_1 &= a\,f(a^5b^3, a^{-1}b) = \tfrac{1}{2}\left[f(a,b) - f(-a,-b)\right],
\end{align*}
which appear as Entry~30(ii)--(iii) in \cite[Part~III, p.~46]{BerndtIII}.
Similarly, at $m=3$, he obtained the cubic transformation (Entry~7, \cite[Part~IV, p.~144]{BerndtIV})
\[
f(\omega a, \omega b) = \omega f(a,b) + (1-\omega)\,f(a^6b^3, a^3b^6),
\]
and at $m=4$ the quartic identity (Entry~9, \cite[Part~IV, p.~146]{BerndtIV})
\[
f(i a,i b) = \tfrac{1}{2}(1+i)f(a,b) + \tfrac{1}{2}(1-i)f(-a,-b).
\]
These transformation formulas reveal deep modular properties of theta functions and are closely related to the transformation behavior of Jacobi's theta functions under the modular group.
The main goal of this paper is to generalize these transformations into a single formula.
\section{Main Results}
\begin{thm}[Generalized Theta Function Transformation]\label{thm:main_transform}
Let $\zeta$ be a primitive $m$-th root of unity. The transformation of $f(a, b)$ is given by the identity:
$$ f(\zeta a, \zeta b) = \sum_{k=0}^{m-1} \zeta^{k^2} a^{k(k+1)/2} b^{k(k-1)/2} f(A_m (ab)^{mk}, B_m (ab)^{-mk}) $$
where
$$ A_m = a^{m(m+1)/2} b^{m(m-1)/2} \quad \text{and} \quad B_m = a^{m(m-1)/2} b^{m(m+1)/2} $$
\end{thm}
\begin{proof}
We begin by applying the definition of $f(a, b)$ to the arguments $(\zeta a, \zeta b)$:
$$ f(\zeta a, \zeta b) = \sum_{n=-\infty}^{\infty} (\zeta a)^{n(n+1)/2} (\zeta b)^{n(n-1)/2} $$
The exponent of $\zeta$ simplifies to $\frac{n(n+1)}{2} + \frac{n(n-1)}{2} = n^2$. Thus,
$$ f(\zeta a, \zeta b) = \sum_{n=-\infty}^{\infty} \zeta^{n^2} a^{n(n+1)/2} b^{n(n-1)/2} $$
Next, we dissect the sum based on the residue of $n$ modulo $m$. Since $n \equiv k \pmod m$ implies $\zeta^{n^2} = \zeta^{k^2}$, we can write:
$$ f(\zeta a, \zeta b) = \sum_{k=0}^{m-1} \zeta^{k^2} \left( \sum_{n \equiv k (\text{mod } m)} a^{n(n+1)/2} b^{n(n-1)/2} \right) $$
Let $S^{(m)}_k(a, b)$ be the inner sum. We analyze $S^{(m)}_k(a, b)$ by substituting $n = mj+k$:
$$ S^{(m)}_k(a, b) = \sum_{j=-\infty}^{\infty} a^{(mj+k)(mj+k+1)/2} b^{(mj+k)(mj+k-1)/2} $$
We factor out the term corresponding to $j=0$ (i.e., $n=k$):
$$ S^{(m)}_k(a, b) = a^{k(k+1)/2} b^{k(k-1)/2} \sum_{j=-\infty}^{\infty} a^{\frac{(mj+k)(mj+k+1)}{2} - \frac{k(k+1)}{2}} b^{\frac{(mj+k)(mj+k-1)}{2} - \frac{k(k-1)}{2}} $$
Simplifying the exponents inside the sum gives:
$$ S_k(a, b) = a^{k(k+1)/2} b^{k(k-1)/2} \sum_{j=-\infty}^{\infty} a^{\frac{m^2j^2 + 2mkj + mj}{2}} b^{\frac{m^2j^2 + 2mkj - mj}{2}} $$
The inner sum is a theta function $f(A_k, B_k)$. By matching the exponents to the definition $f(A_k, B_k) = \sum A_k^{j(j+1)/2} B_k^{j(j-1)/2}$, we solve for $A_k$ and $B_k$:
$$ A_k = a^{(m^2+2mk+m)/2} b^{(m^2+2mk-m)/2} = A_m (ab)^{mk} $$
$$ B_k = a^{(m^2-2mk-m)/2} b^{(m^2-2mk+m)/2} = B_m (ab)^{-mk} $$
Thus, we have shown:
$$ S_k(a, b) = a^{k(k+1)/2} b^{k(k-1)/2} f(A_m (ab)^{mk}, B_m (ab)^{-mk}) $$
Substituting this expression for $S_k(a, b)$ back into our equation for $f(\zeta a, \zeta b)$ completes the proof.
\end{proof}
\begin{cor}[book III,pg.46,Entry 30(ii) and (iii)]
The transformation at $m=2$ allows for the even ($S_0$) and odd ($S_1$) components of the dissection to be isolated:
\begin{align*}
    S_0 &= f(a^3b, ab^3) = \frac{1}{2} \left[ f(a, b) + f(-a, -b) \right] \\
    S_1 &= a f(a^5b^3, a^{-1}b) = \frac{1}{2} \left[ f(a, b) - f(-a, -b) \right]
\end{align*}
\end{cor}
\begin{proof}
From the proof of Theorem \ref{thm:main_transform}, the $m=2$ transformation is given by $f(-a, -b) = S_0 + (-1)^1 S_1 = S_0 - S_1$.
By definition, the $m=2$ dissection is $f(a, b) = S_0 + S_1$.
We now have a simple system of two linear equations:
\begin{align*}
    f(a, b) &= S_0 + S_1 \\
    f(-a, -b) &= S_0 - S_1
\end{align*}
Adding the two equations yields $f(a, b) + f(-a, -b) = 2S_0$, which proves the first identity.
Subtracting the second equation from the first yields $f(a, b) - f(-a, -b) = 2S_1$, which proves the second identity.
\end{proof}
\begin{rem}
 In the above corollary, if you put \(a=b=q\), you get the results in book III, pg. 40, Entry 25(i) and (ii)
\end{rem}
\begin{cor}[book IV,pg.144,Entry 7]
The identity at $m=3$ can be simplified to the compact form:
$$ f(\omega a, \omega b) = \omega f(a, b) + (1 - \omega) f(a^6b^3, a^3b^6) $$
\end{cor}
\begin{proof}
From the proof of Theorem \ref{thm:main_transform}, we know that $f(\omega a, \omega b) = S_0 + \omega S_1 + \omega S_2$, where $S_k$ is the $k$-th term of the dissection.
The three terms in the sum above are precisely $S_0$, $\omega S_1$, and $\omega S_2$.
$$ f(\omega a, \omega b) = S_0 + \omega(S_1 + S_2) $$
By definition, $f(a, b) = S_0 + S_1 + S_2$, which implies $S_1 + S_2 = f(a, b) - S_0$.
Substituting this into the transformation equation:
\begin{align*}
    f(\omega a, \omega b) &= S_0 + \omega(f(a, b) - S_0) \\
    &= S_0 + \omega f(a, b) - \omega S_0 \\
    &= \omega f(a, b) + (1 - \omega) S_0
\end{align*}
Finally, substituting $S_0 = f(A_3, B_3) = f(a^6b^3, a^3b^6)$ gives the result. This simplification is possible because the coefficients for $k=1$ and $k=2$ are identical.
\end{proof}
\begin{cor}[book IV, pg.~146, Entry~9): The $m=4$ transformation and Ramanujan's compact form]
For \(m=4\) and the primitive fourth root of unity \(\zeta = i\),
the transformation of the theta function \(f(a,b)\) satisfies the following equivalent forms:
\\
\textbf{(a) Explicit four-term dissection form:}
\[
\begin{aligned}
f(i a, i b)
 \;=\;&
 f(a^{10} b^{6},\, a^{6} b^{10})
 + a^{3} b\; f(a^{18} b^{14},\, a^{14} b^{18}) \\[3pt]
 &\quad +\, i\!\left[
   a\, f(a^{14} b^{10},\, a^{10} b^{14})
   + a^{6} b^{3}\, f(a^{22} b^{18},\, a^{18} b^{22})
 \right].
\end{aligned}
\]
\textbf{(b) Ramanujan's compact form:}
\[
f(i a,i b)
=\tfrac{1}{2}(1+i)\,f(a,b)\;+\;\tfrac{1}{2}(1-i)\,f(-a,-b).
\]
\end{cor}
\begin{proof}
From Theorem~\ref{thm:main_transform}, for $m=4$ we have
\[
f(i a, i b)
 = \sum_{k=0}^{3} i^{k^2}
 a^{\frac{k(k+1)}{2}} b^{\frac{k(k-1)}{2}}
 f(A_4 (ab)^{4k}, B_4 (ab)^{-4k}),
\]
where $A_4 = a^{10} b^{6}$ and $B_4 = a^{6} b^{10}$.
Evaluating for $k=0,1,2,3$ gives:
\[
\begin{aligned}
k=0&:\quad f(A_4, B_4),\\
k=1&:\quad i\, a\, f(A_4(ab)^4,\, B_4(ab)^{-4})
   = i\, a\, f(a^{14}b^{10},\, a^{10}b^{14}),\\
k=2&:\quad a^{3}b\, f(A_4(ab)^8,\, B_4(ab)^{-8})
   = a^{3}b\, f(a^{18}b^{14},\, a^{14}b^{18}),\\
k=3&:\quad i\, a^{6}b^{3}\, f(A_4(ab)^{12},\, B_4(ab)^{-12})
   = i\, a^{6}b^{3}\, f(a^{22}b^{18},\, a^{18}b^{22}).
\end{aligned}
\]
Collecting the real and imaginary parts gives form (a).
To obtain the compact form (b), recall
\[
f(a,b)=S_0+S_1+S_2+S_3, \quad f(-a,-b)=S_0-S_1+S_2-S_3,
\]
and
\[
f(i a,i b)=(S_0+S_2) + i(S_1+S_3).
\]
Eliminating $S_k$ gives
\[
S_0+S_2=\tfrac{1}{2}(f(a,b)+f(-a,-b)),\quad
S_1+S_3=\tfrac{1}{2}(f(a,b)-f(-a,-b)).
\]
Substituting back yields
\[
\begin{aligned}
f(i a,i b)
&= \tfrac{1}{2}(f(a,b)+f(-a,-b))
   + \tfrac{i}{2}(f(a,b)-f(-a,-b))\\[4pt]
&= \tfrac{1}{2}(1+i)\,f(a,b)
   + \tfrac{1}{2}(1-i)\,f(-a,-b),
\end{aligned}
\]
which is the compact form (b), equivalent to (a).
\end{proof}
\begin{rem}[Equating real and imaginary parts]
Comparing the two equivalent forms (a) and (b), here we use corollary 2.5 also,
\[
\begin{aligned}
\Re\!\bigl(f(i a,i b)\bigr)
&=\tfrac{1}{2}\!\left[f(a,b)+f(-a,-b)\right]
 = f(a^{10} b^{6},\, a^{6} b^{10}) + a^{3} b\; f(a^{18} b^{14},\, a^{14} b^{18})=f(a^3b, ab^3),\\[6pt]
\Im\!\bigl(f(i a,i b)\bigr)
&=\tfrac{1}{2}\!\left[f(a,b)-f(-a,-b)\right]
 = a\, f(a^{14} b^{10},\, a^{10} b^{14}) + a^{6} b^{3}\, f(a^{22} b^{18},\, a^{18} b^{22})=a f(a^5b^3, a^{-1}b).
\end{aligned}
\]
Hence, the \emph{real part} corresponds to the even dissection components $(S_0,S_2)$, while the \emph{imaginary part} corresponds to the odd components $(S_1,S_3)$.
Setting $a=b=q$ gives
\[
\begin{aligned}
\Re[f(iq,iq)] &= f(q^{16}, q^{16}) + q^4 f(q^{32}, q^{32}) = \tfrac{1}{2}\!\left[f(q,q)+f(-q,-q)\right],\\[4pt]
\Im[f(iq,iq)] &= q f(q^{24}, q^{24}) + q^{10} f(q^{40}, q^{40}) = \tfrac{1}{2}\!\left[f(q,q)-f(-q,-q)\right],
\end{aligned}
\]
\end{rem}
\section{acknowledgment}
The author thanks Professor A.K.Shukla for his encouragement and constant support.

\end{document}